\documentclass[11pt,twoside]{amsart}
\usepackage{amscd}

\theoremstyle{plain}
\newtheorem{thm}{Theorem}[section]

\newtheorem{lem}[thm]{Lemma}

\newtheorem{prop}[thm]{Proposition}

\theoremstyle{definition}

\newtheorem*{ackno}{Acknowledgments} 

\theoremstyle{remark}

\author{Vitaly Vologodsky} 
\address{Department of Mathematics\\
University of Georgia\\
Athens, GA 30602}
\email{vologods@math.uga.edu}

\begin{document}
\bibliographystyle{amsalpha}
\title[Prym Map]%
{Locus of Indeterminacy of the Prym Map} 
\date{March 29, 2001}

\begin{abstract}
We provide an easy characterization for 
the locus of indeterminacy of the Prym map 
in terms of the dual graphs of stable curves.
As a corollary, we show that the closure of the 
Fridman-Smith locus coincides with the locus 
of indeterminacy of the Prym map.
\end{abstract}

\maketitle

\setcounter{section}{-1}

\section{Introduction}
\label{sec:Introduction}

This paper is intended to be an appendix to the \cite{ABH}
and here we are using notations and results of \cite{ABH} .

The extended Prym map is a rational map from 
the moduli space $\bar {R_g}$
of stable curves $C$ of genus $2g-1$ with an involution $i$, 
where the only base points of $i$ are nodes and the involution does
not exchange branches at that node,
to the space 
$\bar A^{vor}_{g-1}$, the
toroidal compactification of $A_{g-1}$ for the 2nd Voronoi fan.
On the dense open subset in $\bar R_g $
corresponding to smooth curves it is given by associating to a pair
$(C, i)$ its Prym variety $P(C,i)$.
Locus of indeterminacy of this map was described
in \cite{ABH} by a condition (*)
using the combinatorial data of the dual graph of a curve with an involution.
This paper gives an easier description for the indeterminacy locus.

\begin{thm}\label{thm:intro}
A curve with an involution is in the indeterminacy locus
of the extended Prym map if and only if  
the curve is a degeneration of a Friedman-Smith
example with the number of edges 4 or greater.
\end{thm}

Recall that the Friedman-Smith example with $2n$ edges is a curve 
$(C,i)\in \bar {R_g}$, such that the curve $C= C_1 \cup C_2$
is the union of two irreducible components, both invariant
under the involution, intersecting in $2n$ points, so
that the involution is base point free and interchanges these
$2n$ nodes.

For $n\ge 2$ the closure of the Friedman-Smith locus with $2n$ edges
contains $[\frac {g-n+1}2]+1$ irreducible components. Curve $C'=C/i$ is the
union $C'=C_1'\cup C_2'$ of two irreducible curves of genuses
$g_1'$ and $g_2'$, where $g_1'+g_2'=g-n+1$.
Irreducible components of the Friedman-Smith 
locus correspond to cases when $g_1'$ and $g_2'$
are $(0, g-n+1), (1, g-n), \dots ,([\frac {g-n+1}2],[\frac {g-n+2}2])$.

Note that in our space $\bar {R_g}$ there are no stable curves 
having double points with the 2 branches exchanged 
by the involution (called the nodes of type (2) in \cite{ABH}).

\begin{ackno}
The author would like to thank Professor V.Alexeev 
for numerous helpful discussions and Professor K.Hulek
for valuable comments on the preliminary version of this paper.
\end{ackno}

\medskip

\section{Preliminary results}
\label{sec:Preliminary results}

We will call a vertex of the dual graph a {\it bold vertex} if 
the corresponding component of the curve is fixed by the involution
and we will call an edge a {\it bold edge} if the
corresponding node is a node of type (1), see
\cite{ABH} for definition.
We will call all other vertices and edges {\it ordinary}. 
Note that if we have a bold edge then its beginning and the end
are bold vertices. Therefore, the union of bold vertices and edges 
forms a subgraph of the dual graph. We will denote this subgraph
by $B(\Gamma)$.

Proof of one direction for the Theorem \ref{thm:intro} is immediate.

\begin{prop}\label{prop:easyway}
If the curve is a degeneration of a Friedman-Smith example with
at least 4 edges, then the curve is in the indeterminacy locus 
of the Prym map.
\end{prop}
\begin{proof}
Prym map is not defined for the Friedman-Smith examples, 
see \cite{FS} or \cite{ABH} Section 4.2.3.
The indeterminacy locus for every rational map is closed,
therefore degenerations of the Friedman-Smith examples 
are also in the indeterminacy locus.
\end{proof}

Let us consider the dual graphs which are not degenerations of
Friedman-Smith examples with the number of edges at least 4.

\begin{lem}\label{lem:4edg}
If the dual graph $\Gamma$
of a stable curve $C$ with an involution contains
two disjoint connected equivariant subgraphs 
$\Gamma _1$  and $\Gamma _2$
which are connected 
by at least 4 ordinary edges (so that the beginnings of these edges are 
in one subgraph and the ends are in the other), 
but are not connected by any bold path
then this curve is a degeneration of a Friedman-Smith example 
with at least 4 edges.
\end{lem}

\begin{proof}
We will construct two disjoint  equivariant subgraphs 
$\Gamma '_1$ and $\Gamma '_2$
such that together they contain all vertices of the graph $\Gamma$ 
and every edge of $\Gamma$ either belongs to one of the subgraphs
$\Gamma '_1$ or $\Gamma '_2$
or is an ordinary edge which begins at one subgraph and ends at the 
other. 

To find these $\Gamma '_1$ and $\Gamma '_2$ we will take 
$\Gamma _1$  and $\Gamma _2$ and start adding to them 
vertices and edges until we will get required $\Gamma '_1$ and $\Gamma '_2$.

First let us consider connected components of the bold 
subgraph $B(\Gamma)$.
Each bold component cannot intersect with both subgraphs
$\Gamma _1$ and $\Gamma _2$, otherwise we have a bold path
connecting these subgraphs. We will change each $\Gamma _i$ by
the union of $\Gamma _i$ and all components of $B(\Gamma )$
which has non empty intersection with $\Gamma _i$.
It is easy to see that the new $\Gamma _1$ and $\Gamma _2$
still satisfy all the conditions from the statement of the Lemma.
In addition, the compliment 
$\Gamma \backslash (\Gamma _1 \cup \Gamma _2 )$
is attached to $\Gamma _1 \cup \Gamma _2 $
through the ordinary edges only.

Next let us consider the following topological space 
$\Gamma \backslash (\Gamma _1 \cup \Gamma _2 )$.
Connected components of this space 
can be of three different types:

1) components attached to $\Gamma _1$ only,

2) components attached to $\Gamma _2$ only,

3) components attached to both $\Gamma _1$ and $\Gamma _2$.

We will change $\Gamma _1$ by the union 
$\Gamma _1 \cup (\text{all components of type 1})$
and $\Gamma _2$ by
$\Gamma _2 \cup (\text{all components of type 2})$.
Again, $\Gamma _1$ and $\Gamma _2$
still satisfy all the conditions from the statement of the Lemma.

Now, to get desired $\Gamma '_1$ and $\Gamma '_2$,
we can take $\Gamma '_1$ to be equal $\Gamma _1$.
For vertices of $\Gamma '_2$ we will take all vertices of 
$\Gamma$ which are not from subgraph $\Gamma '_1$
and for edges of $\Gamma '_2$ we will take all edges 
whose beginning and ends are vertices of $\Gamma '_2$.

Note that the number of edges connecting $\Gamma '_1$ and $\Gamma '_2$
is at least the same as the number of edges connecting original
$\Gamma _1$ and $\Gamma _2$, which is not less then 4.

Consider decomposition of the curve $C=C_1 \cup C_2$  
where each $C_i$ correspond to subgraph $\Gamma '_i$.
Each part $C_i$ is invariant with respect to involution
and can be obtained as a degeneration of smooth curves with involution.
Then the whole curve $C$ is a degeneration of the Friedman-Smith example 
and as we have noticed the number of nodes is at least 4.
\end{proof}

\bigskip

\section{Condition (*)}
\label{sec:*}

Let us fix an orientation on the graph $\Gamma$ which is compatible 
with the involution $i$.
Following \cite{ABH} edges $e_j$ give coordinate functions on
$C_1(\Gamma , \mathbb {R} )$. We have three possibilities:

1. $z_j$ is identically zero on $X^{-}$.

2. Image $z_j(X^{-})=\mathbb {Z}$.

3. Image $z_j(X^{-})=\frac 12 \mathbb Z$. 

\medskip
${}$ From \cite{ABH} we have the following combinatorial condition
\medskip

(*) The linear functions $m_j z_j$ define 
a dicing of the lattice $X^{-}$.
\medskip

In this dicing we are not using functions of type 1,
$m_j=1$ for the functions $z_j$ of type 2 
and $m_j=2$ for the functions $z_j$ of type 3. 

Another characterization of possibilities 2 and 3 for image $z_j(X^-)$ 
is as follows. 
Image $z_j(X^{-})=\mathbb {Z}$ if for every simple oriented cycle
$\omega \in H^1(\Gamma , \mathbb Z )$ in the graph $\Gamma $
$\operatorname{mult}_{e_j} \omega =1$
implies $\operatorname{mult}_{e_{i(j)}} \omega =-1$.
Image $z_j(X^{-})=\frac 12\mathbb {Z}$ if 
there exists a simple oriented cycle 
$\omega \in H^1(\Gamma , \mathbb Z )$ with
$\operatorname{mult}_{e_j} \omega =1$
but $\operatorname{mult}_{e_{i(j)}} \omega =0$.

Let $d$ be the dimension of the space $X^-\otimes \mathbb R$.
Then condition (*) is equivalent to the following:
\medskip

If the intersection of $d$ hyperplanes 
$\{m_{j_k} z_{j_k} =n_{j_k}, n_{j_k}\in \mathbb Z, k=1,\dots ,d\}$
is 0-dimensional (i.e. a point),
then it is in the lattice $X^-$.
\medskip

It is enough to check this condition for the sets of 
$\{n_{j_k}, k=1,\dots ,d\}$ where all except one $n_{j_k}$ are 0
and the remaining $n_{j_k}$ is 1. 
If the intersections of hyperplanes corresponding to these sets 
$\{n_{j_k}, k=1,\dots ,d\}$ are in the lattice $X^-$ then,
since all other $\{n_{j_k}, k=1,\dots ,d\}$
are integer linear combinations of those sets, their
corresponding intersections of the hyperplanes are 
integer linear combinations
of elements from $X^-$, and therefore belong to $X^-$.

We can reprove Proposition \ref{prop:easyway} without using 
the geometrical meaning of condition (*).

\begin{proof}
Let curve $C$ be a degeneration of the 
Friedman-Smith example with $2n$ edges with $n \ge 2$. 
Thus curve $C$ is a union $C_1\cup C_2$ where 
each $C_j$ is a degeneration of a family of smooth curves 
and intersection $C_1\cap C_2$ consists of $n$ pairs of exchanged nodes.
Then in the dual graph for $C$ there are two equivariant subgraphs
$\Gamma_1$ and $\Gamma_2$ corresponding to $C_1$ and $C_2$
and these subgraphs
$\Gamma_1$ and $\Gamma_2$ are connected by the $n$ pairs of 
exchanging edges.
We will denote them $e_1, e_1 '=i(e_1), \dots , e_n, e_n '=i(e_n)$.
It is clear that these edges are of type 3. In definition of
condition (*) all $m_j$ corresponding to 
$e_1, e_1 ', \dots , e_n, e_n '$ equal 2.

Now if we have an arbitrary simple cycle $\omega$ then either
it does not contain any of
the edges $e_1, e_1 ', \dots , e_n, e_n '$
or it passes through even number of them.
With counting multiplicities the sum $\omega -i(\omega )$ 
passes through $e_1, \dots , e_n$ the same even number of times.
All the elements of the lattice $X^-$ are of the form
$\frac {1-i}2 (\sum\limits_l k_l\omega _l)$ where 
$k_l \in \mathbb Z$ and $\omega _l$ are
simple oriented cycles.
Therefore for any element of $X^-$  we have
$\sum\limits_{m=1}^n
{\operatorname{mult}_{e_m}}
\big(\frac {1-i}2 (\sum\limits_l k_l\omega _l)\big)$
is an integer.
We can choose some $d-n$  
edges $e_{j_{n+1}}, \dots ,e_{j_d}$
with the functions $z_{j_{n+1}}, \dots ,z_{j_d}$
so that the intersection of $d$ hyperplanes
$\{2z_1 =1, 2z_2=0,\dots , 2z_n=0, m_{j_{n+1}}z_{j_{n+1}}=0, 
\dots, m_{j_d}z_{j_d}=0\}$
is 0-dimensional. But the point of intersection $x$ cannot be 
an element of $X^-$ since
$\sum\limits_{m=1}^n
{\operatorname{mult}_{e_m}}(x)=\frac 12$ is not an integer.
Condition (*) fails.
\end{proof}

To finish the proof of Theorem \ref{thm:intro}, let us prove the other direction 
in its statement.

Let us compare functions $z_j$ and $z_{i(j)}$ when restricted to $X^-$.
Elements of $X^-$ are of the form
$\frac {1-i}2 (\sum\limits_l k_l\omega _l)$ where 
$k_l \in \mathbb Z$ and $\omega _l$ are
simple oriented cycles. 
If $\omega _l = \sum\limits_{m=1}^{\operatorname{length}_{\omega _l}}
    e_{l_m}$
then
$$
\frac {1-i}2 \big[\sum\limits_l k_l\omega _l\big]=
\frac 12 \big[\sum\limits_l k_l \big(
\sum\limits_{m=1}^{\operatorname{length}_{\omega _l}}e_{l_m} \big)
-\sum\limits_l k_l \big(
\sum\limits_{m=1}^{\operatorname{length}_{\omega _l}}i(e_{l_m}) \big)
\big].
$$
Since the orientation of the graph $\Gamma$ 
is compatible with the involution $i$
we have $i(e_{l_m})=+e_{i(l_m)}$ and thus $z_j=-z_{i(j)}$
when restricted on $X^-$.

Intersection of $d$ hyperplanes 
$\{m_{j_k} z_{j_k} =n_{j_k}, n_{j_k}\in \mathbb Z, k=1,\dots ,d\}$
is 0-dimensional if and only if the intersection of
$\{m_{j_k} z_{j_k} =0, k=1,\dots ,d\}$ is 
$0\in X^- \otimes R$. The later is equivalent to 
$\{m_{j_k} z_{j_k} =0, m_{j_k} z_{i(j_k)} =0, k=1,\dots ,d\}$ 
having intersection $0\in X^- \otimes R$.
Note that $z_j(\omega )=z_{i(j)}(\omega )=0$ 
for an element $\omega \in X^-$ is equivalent to the fact
that $\omega $ can be represented as 
$\frac {1-i}2 \big( \sum\limits_l k_l\omega _l\big)$ where all
$\omega _l$ are simple cycles in $\Gamma$ not passing through $e_j$ or
$e_{i(j)}$.
Therefore intersection of
$\{m_{j_k} z_{j_k} =0, m_{j_k} z_{i(j_k)} =0, k=1,\dots ,n\}$
consist of the elements of $X^-$ of the form
$\frac {1-i}2 \big( \sum\limits_l \omega _l\big)$ with
$\omega _l \in H^1(\Gamma\backslash 
(e_{j_1}\cup e_{i(j_1)} \cup\dots\cup e_{j_n}\cup e_{i(j_n)}),
\mathbb Z )$.
Finally, intersection of $d$ hyperplanes 
$\{m_{j_k} z_{j_k} =n_{j_k}, n_{j_k}\in \mathbb Z, k=1,\dots ,d\}$
is 0-dimensional if and only if for the graph
$\Gamma '=\Gamma\backslash
(e_{j_1}\cup e_{i(j_1)} \cup\dots\cup e_{j_n}\cup e_{i(j_n)})$
all elements $\omega \in H^1(\Gamma ',\mathbb Z)$
are invariant with respect to $i$.

In the case of 0-dimensional intersection let us look closely 
at the subgraph $\Gamma '$. 
It may have more than one connected component.
But it is impossible that two of these components are exchanged
by the involution. Indeed, for connected curve with involution
dimension of the corresponding $X^-$ can be computed as $n_e-c_e$
where $2n_e$ is the number of exchanged under the involution 
nodes (nodes of type (3)) 
(or the number of ordinary edges in the dual graph)
and $2c_e$ is the number of exchanged components
(or the number of ordinary vertices).
For the graph $\Gamma $ we have $n_e-c_e=d$. Subgraph
$\Gamma '$ was obtained by erasing $2d$ ordinary edges, thus for
$\Gamma '$ we have $n_e-c_e=0$.
Now $X^-$ for $\Gamma '$ is 0, so for every 
invariant under the involution connected component
of $\Gamma '$ we must have $n_e-c_e=0$.
Components of $\Gamma '$ exchanged by $i$ consist of
ordinary vertices and ordinary edges only.
To have $X^-=0$ for $\Gamma '$
all exchanged components should have no cycles.
Then they are trees and for each of them $n_e-c_e=-1$.
If the exchanged components present in $\Gamma '$
then for $\Gamma '$ $n_e-c_e$ would be negative.
All connected components of $\Gamma '$ are fixed under $i$.

Now consider the point of intersection for the hyperplanes 
$\{m_{j_1}z_{j_1}=1, m_{j_k}z_{j_k}=0, k=2,\dots ,d\}$.

If the edge $e_{j_1}$ has the beginning $v_{j_1}$ 
and the end $v_{j_1}'$ at the same
connected component of $\Gamma '$ then there is a path $t$
inside of $\Gamma '$ which starts at 
$v_{j_1}'$ and ends at $v_{j_1}$.
Cycle $\omega =t+e_{j_1}$ has 
$\operatorname{mult}_{e_{j_1}} \omega =1$
but $\operatorname{mult}_{e_{i(j_1)}} \omega =0$
which means that $m_{j_1}=2$.
Now the element $x^-=\frac {1-i}2 (\omega ) \in X^-$ satisfies 
$\{m_{j_1}z_{j_1}(x^-)=1, m_{j_k}z_{j_k}(x^-)=0, k=2,\dots ,d\}$.
The point of intersection of the hyperplanes is in $X^-$.

If the edge $e_{j_1}$ has the beginning $v_{j_1}$
and the end $v_{j_1}'$ in different components of $\Gamma '$
then the edges $e_{j_1}$ and $e_{i(j_1)}$ 
start (resp. end) at the same component. 
There is a path $t$ inside of $\Gamma '$ which starts at
$v_{j_1}$ and ends at $i(v_{j_1})$. 
Similarly, there is a path
$t'$ inside of $\Gamma '$ which starts at
$v_{j_1}'$ and ends at $i(v_{j_1}')$.
$\omega =e_{j_1}+t'-e_{i(j_1)}-t$ is a cycle in $\Gamma$.
Element $x^-=\frac {1-i}2 (\omega ) \in X^-$
satisfies 
$\{z_{j_1}(x^-)=1, z_{j_k}(x^-)=0, k=2,\dots ,d\}$.
To see that $x^-$ is
the point of intersection for the hyperplanes 
$\{m_{j_1}z_{j_1}=1, m_{j_k}z_{j_k}=0, k=2,\dots ,d\}$
we need to show that $m_{j_1}=1$.

We will show that if 
there exists a simple oriented cycle 
$\omega \in H^1(\Gamma , \mathbb Z )$ with
$\operatorname{mult}_{e_{j_1}} \omega =1$
but $\operatorname{mult}_{i(e_{j_1})} \omega =0$ then 
$\Gamma$ is a degeneration of Friedman-Smith example with 
at least 4 edges.
Cycle $\omega $ goes along some of the edges $e_{j_k}$ and
$e_{i(j_k)}$ and passes through some connected components
$C_1,\dots C_l$ of the subgraph $\Gamma '$. We can assume that the
edge $e_{j_1}$ starts at $C_1$ and ends at $C_2$.
Cycle $\omega $ is simple so it cannot reenter $C_1$
along $e_{j_1}$. $\operatorname{mult}_{i(e_{j_1})} \omega =0$, so
$\omega $ cannot reenter $C_1$ along $i(e_{j_1})$.
There is another edge $e'$ through which $\omega $ reenters $C_1$.
Let $\Gamma _1=C_1$, $\Gamma _2'=(C_2\cup\dots\cup C_l)\cup 
(\omega\backslash(C_1\cup e_{j_1}\cup e'))$ and 
$\Gamma _2=\Gamma_2'\cup i(\Gamma_2')$.
Subgraphs $\Gamma _1$ and $\Gamma _2$ are connected by at least
4 ordinary edges: $e_{j_k}$, $e_{i(j_k})$, $e'$ and $i(e')$.
Therefore, they satisfy the assumptions of the Lemma \ref{lem:4edg} 
and graph $\Gamma$ is a degeneration of Friedman-Smith example with 
at least 4 edges.

This proves Theorem \ref{thm:intro}.

\bigskip

\section{Condition (**)}
\label{section:**}

In this section we will discuss condition (**) of \cite{ABH}
which is stronger than 
condition (*).
When condition (*) holds but (**) fails we have a Prym variety 
corresponding to the double cover of stable curves but this Prym 
cannot be canonically embedded into the Jacobian of the covering curve. 

${}$From \cite{ABH} we have the following definition of
combinatorial condition (**)
\medskip

(**) $[\Delta ^- ]$ is a dicing with respect to the lattice
$2X^-$.
\medskip

In the notations of the previous
section this condition (**) is equivalent to the following:\medskip

If the intersection of $d$ hyperplanes 
$\{z_{j_k} =n_{j_k}, n_{j_k}\in \mathbb Z, k=1,\dots ,d\}$
is 0-dimencional (i.e. a point),
then it is in the lattice $2X^-$.
\medskip

The first condition of the following theorem gives an
easier description for (**).

\begin{thm}\label{thm:**}
The following three conditions are equivalent

(i) Curve is not a degeneration of a Friedman-Smith
example with the number of edges 2 or greater.

(ii) Condition (*) holds and 
there are no functions $z_j$ of type (2).

(iii) Condition (**) holds.

\end{thm}
\begin{proof}
{\it (i) implies (ii).} If (*) does not hold then
(**) does not hold either.
We only need to show that if (*) holds but the graph has 
an edge $e_j$ with the corresponding function
$z_j$ of type (2) then the graph is a degeneration of 
a Friedman-Smith example with 2 edges.

Let $v_j$ be the beginning of $e_j$ and $v_j'$ be the
end of $e_j$. Let us consider connected components of
$\Gamma_j'=\Gamma \backslash (e_j \cup e_{i(j)} )$. 
The number of components is at most 3.
If there is only one component in $\Gamma_j'$
then there is a path $t$ in $\Gamma_j'$ connecting
$v_j$ and $v_j'$. Thus for the simple cycle 
$\omega = e_j -t$ me have 
$\operatorname{mult}_{e_{j}} \omega =1$ and
$\operatorname{mult}_{e_{i(j)}} \omega =0$
which is impossible since $e_j$ is an edge of type 2.
If there were three components in $\Gamma_j'$
then edges $e_j$ and $e_{i(j)}$ would have been of type 1.
The only possibility is when $\Gamma_j'$ has two components 
$\Gamma_{j1}'$ and $\Gamma_{j2}'$. Since for every simple cycle 
$\omega$ we have
$\operatorname{mult}_{e_{j}} \omega =
-\operatorname{mult}_{e_{i(j)}} \omega$
then beginnings $v_j$ and $v_j'$ belong to the same component
of $\Gamma_j'$, say
$\Gamma_{j1}'$. Therefore $\Gamma_{j1}'$ is invariant with respect
of involution $i$ and so is $\Gamma_{i2}'$.
Now consider decomposition of the curve $C=C_1 \cup C_2$  
where each $C_k$ correspond to subgraph $\Gamma '_{jk}$.
Each part $C_k$ is invariant with respect to involution
and can be obtained as a degeneration of smooth curves with involution.
Then the whole curve $C$ is a degeneration of the Friedman-Smith example 
with the number of edges 2.

{\it (ii) implies (i).} Assume that (i) fails.
The curve $C$ is a degeneration of a 
Friedman-Smith example with 2 edges. 
Thus curve $C$ is a union $C_1\cup C_2$ where 
each $C_j$ is a degeneration of a smooth curve 
and intersection $C_1\cap C_2$ consists of two exchanged nodes.
Then in the dual graph for $C$ there are two equivariant subgraphs
$\Gamma_1$ and $\Gamma_2$ corresponding to $C_1$ and $C_2$
and these subgraphs
$\Gamma_1$ and $\Gamma_2$ are connected by the pair of 
exchanging edges.
It is clear that these edges are of type 2.

{\it (ii) implies (iii).} Since there are no edges of type 2
then in definition of condition (*) all $m_j$ equal 2.
Condition (*) holds, so
if the intersection of $d$ hyperplanes 
$\{2 z_{j_k} =n_{j_k}, n_{j_k}\in \mathbb Z, k=1,\dots ,d\}$
is 0-dimensional,
then it is in the lattice $X^-$. 
Thus intersection of 
$\{2 z_{j_k} =2n_{j_k}, n_{j_k}\in \mathbb Z, k=1,\dots ,d\}$
is in $2X^-$, but this is the same point
as intersection of
$\{z_{j_k} =n_{j_k}, n_{j_k}\in \mathbb Z, k=1,\dots ,d\}$.
Therefore condition (**) holds.

{\it (iii) implies (ii).} Let us take an arbitrary edge $e_j$
and the corresponding function $z_j$. If $z_j$ is nontrivial
(i.e. not of type 1) then we can choose some $d-1$  
edges $e_{j_2}, \dots ,e_{j_d}$
with the functions $z_{j_2}, \dots ,z_{j_d}$
so that the intersection of $d$ hyperplanes
$\{z_{j} =1, z_{j_2}=0,\dots , z_{j_d}=0\}$
is 0-dimensional and by condition (**) is
some point $2x^- \in 2X^-$. Now $x^- \in X^-$
and $z_j(x^-)=\frac 12$. Thus edge $e_j$
and the corresponding function $z_j$ are of type 3
and not of type 2.

\end{proof}

\end{document}